\documentclass[10pt,a4paper]{article}
\usepackage[utf8]{inputenc}
\usepackage[english]{babel}
\usepackage{amsmath}
\usepackage{amsfonts}
\usepackage{amssymb}
\usepackage{amsthm}
\usepackage{graphicx}
\usepackage{cite}
\usepackage{algorithmic}
\usepackage{textcomp}
\usepackage{xcolor}
\usepackage{paralist}
\usepackage{hyperref}

\newcommand{\set}[1]{\left\lbrace #1 \right\rbrace}
\newcommand{\norm}[1]{\left\Vert #1 \right\Vert}
\newcommand{\abs}[1]{\left\vert #1 \right\vert}

\newcommand{\simto}{\xrightarrow{\raisebox{-0.7ex}[0ex][0ex]{$\sim$}}} 

\theoremstyle{definition}
\newtheorem{theorem}{Theorem}[section]
\newtheorem{proposition}[theorem]{Proposition}

\newtheorem{example}[theorem]{Example}
\newtheorem{definition}[theorem]{Definition}
\newtheorem{lemma}[theorem]{Lemma}
\newtheorem{definition-proposition}[theorem]{Definition-Proposition}
\newtheorem{remark}[theorem]{Remark}

\newtheorem{assumption}[theorem]{Assumption}

\begin{document}

\title{Towards funnel control on Riemannian manifolds
\thanks{This work was funded by the Deutsche Forschungsgemeinschaft (DFG, German Research Foundation) - Project-ID 544702565.}
}

\author{Jonas Kirchhoff$^\dagger$ and Thomas Berger%
\thanks{J. Kirchhoff and T. Berger are with the Institute of Mathematics, Martin Luther University Halle-Wittenberg, Halle (Saale), Germany
        {\tt\small \{jonas.kirchhoff, thomas.berger\}@mathematik.uni-halle.de}}%
}

\maketitle	
\thispagestyle{empty}
\pagestyle{empty}

\begin{abstract}
Funnel control for time-invariant control-affine systems with relative degree one is transferred to the case of dynamical systems evolving on a Riemannian manifold. It is shown that, under some regularity assumptions, the classical proofs carry over. The regularity assumptions reflect the particularities of the nonlinear theory and degenerate in the classical Euclidean case.
\end{abstract}

\section{Introduction}

In this note, we consider the problem of reference tracking with prescribed transient behaviour for the control system 
\begin{align}\label{eq:toy_problem}
\tfrac{\mathrm{d}}{\mathrm{d}t}x = f(x)+B(x)u
\end{align}
where $x$ evolves on a smooth manifold $M$ and the input $u$ is in a not necessarily trivial vector bundle $E\to M$; $f\in\mathfrak{X}(M)$ is a vector field and $B\in\Gamma(\mathcal{T}M\otimes E^*)$ is a vector bundle morphism. The control objective is to find a suitable state feedback control law so that the vector field of the closed loop is forward complete and the distance of the state to a given reference trajectory is always bounded by a given (time-dependent) bound. Here, the distance to the reference is measured in a given Riemannian metric on $M$. The control law should utilise as little information on $f$ and $B$ as possible.

In case that $M$ is Euclidean (i.e. $M=\mathbb{R}^n$ and the Riemannian distance is the Euclidean metric) and $E$ is trivial ($E = M\times\mathbb{R}^m$), it is well-known that the funnel controller achieves the control objective, cf.~\cite{IlchRyanSang02,IlchRyan08,BergIlchRyan21,BergIlchRyan25}. In the general nonlinear case, however, there are two obstructions to the applicability of the Euclidean funnel controller: the absence of global coordinates, or even local ones containing both trajectory and reference, and the nonlinearity of the metric.

Therefore, we shall make an attempt to generalise the Euclidean funnel controller to Riemannian manifolds. The control-affine system~\eqref{eq:toy_problem} takes the r\^{o}le of the linear prototype of Euclidean funnel control. For simplicity of exposition, we consider the relative-degree-one case only, i.e. we assume that $B$ is a vector-bundle-isomorphism. Without loss of generality, we consider connected manifolds only. Moreover, we shall entirely work in the smooth category, with the understanding that smoothness may (if necessary) be replaced by sufficient regularity so that the objects in the proofs may be defined in whatever strong or weak sense required.

\section{Reminder on differential-geometric notions}

In this section, we recall elementary notions from Riemannian geometry used throughout the note. 

\begin{definition}[cf.~\cite{GreuHalpVans72,Lee18}]
A \textit{vector bundle metric} on the vector bundle $E\to M$ is a symmetric tensor field $g\in\Gamma(E^*\vee E^*)$ so that $\langle g, e\otimes e\rangle> 0$ for all $e\in E\setminus 0_E$; $0_E$ is the zero-section. $\norm{e}_g := \sqrt{\langle g,e\otimes e\rangle}$ is the induced fibrewise norm. A \textit{Riemannian metric} on the manifold $M$ is a vector bundle metric $g$ on the tangent bundle $\mathcal{T}M$. The induced Riemannian distance is the infimal distance along piecewise differentiable paths, i.e.
\begin{small}
\begin{align*}
&d_g:M\times M  \to\mathbb{R},\\
&d_g(x,y)  := \inf\set{\int_{t_0}^{t_1} \norm{\tfrac{\mathrm{d}}{\mathrm{d}t}f(t)}_g\,\mathrm{d}t\left\vert\!\!\begin{array}{l}
t_0\leq t_1\\
f\in\mathcal{C}_{\mathrm{pw}}^1([t_0,t_1],M),\\
f(t_0) = x, f(t_1) = y
\end{array}\!\!\!\right.}.
\end{align*}
\end{small}
\end{definition}

For parallelisable manifolds, e.g. open sets $M\subseteq\mathbb{R}^n$ and $S^1$, Riemannian metrics are positive definite matrix fields.

\begin{example}
\begin{enumerate}[(i)]
\item Let $M\subseteq\mathbb{R}^n$ be an open set. $\mathcal{T}M$ is identified with $M\times\mathbb{R}^n$ by means of the canonical basis. A Riemannian metric on $M$ is given by a smooth matrix field $g\in\mathcal{C}^\infty(M,\mathbb{R}^{n\times n})$ so that $g(x)$ is symmetric and positive definite for all $x\in M$. In general, there is no closed form for the Riemannian distance induced by an arbitary Riemannian metric. In case that $g$ is constant,
\begin{align*}
d_g(x,y) = \sqrt{(x-y)^\top g(x-y)}
\end{align*}
with the standard Euclidean metric given by $g\equiv I_n$.
\item Let $M = S^1$ be the one-dimensional sphere. $\mathcal{T}S^1$ is globally trivialisable, and by $S^1 = \mathbb{R}/2\pi\mathbb{Z}$, a Riemannian metric may be identified with an everywhere positive $2\pi$-periodic function $\lambda$, where $g = \lambda^2$ is the tensor field. The Riemannian distance between two points $[\vartheta]$ and $[\vartheta']$ with $\vartheta\leq\vartheta'\in[0,2\pi]$ is then
\begin{align*}
d_g([\vartheta],[\vartheta'])\! =\! \min\set{\int_{\vartheta}^{\vartheta'}\!\lambda(\alpha)\,\mathrm{d}\alpha,L-\int_{\vartheta}^{\vartheta'}\!\lambda(\alpha)\,\mathrm{d}\alpha}\!,
\end{align*}
where $L = \int_0^{2\pi}\lambda(\alpha)\,\mathrm{d}\alpha$. The standard metric is given by $\lambda\equiv 1$.
\end{enumerate}
\end{example}

When we examine Euclidean funnel control, then the control law is given by the difference of the trajectory and the reference multiplied with a gain, which depends on the distance of the two. On a Riemannian manifold, the latter is given by the Riemannian distance, and the former is recognised as the direction of the shortest path connecting the state and its reference. Shortest paths in Riemannian geometry are geodesics, and the direction of a (unique) geodesic is the inverse of the Riemannian exponential. The precise definition is given in terms of the Levi-Civita connection.

\begin{definition}[{cf.~\cite{Lee18}}]
A \textit{connection} on $M$ is a map
\begin{align*}
\nabla:\mathfrak{X}(M)\times\mathfrak{X}(M)\to\mathfrak{X}(M),\qquad (X,Y)\mapsto\nabla_XY
\end{align*}
so that
\begin{enumerate}[(i)]
\item $\nabla$ is $\mathcal{C}^\infty(M)$-linear in the first argument
\item for each $x\in\mathfrak{X}(M)$, $\nabla_X$ is a first order differential operator with symbol $X$, i.e.
\begin{align*}
\nabla_X(\lambda Y) = \lambda \nabla_XY+\langle X,\lambda\rangle Y
\end{align*}
for all $Y\in\mathfrak{X}(M)$ and $\lambda\in\mathcal{C}^\infty(M)$.
\end{enumerate}
The \textit{Levi-Civita connection} induced by the Riemannian metric $g$ is the unique connection $\nabla$ satisfying
\begin{enumerate}[(a)]
\item $\langle X,\langle g,Y\otimes Z\rangle\rangle = \langle g,\nabla_X Y\otimes Z+Y\otimes\nabla_XZ\rangle$
\item $\nabla_XY-\nabla_YX = [X,Y]$, where $[\cdot,\cdot]$ denotes the Lie bracket of vector fields
\end{enumerate}
\end{definition}

Geodesics are curves with vanishing acceleration. The latter is defined with respect to the Levi-Civita connection.

\begin{definition-proposition}[cf.~\cite{Lee18}]
Let $(M,g)$ be a Riemannian manifold with Levi-Civita connection $\nabla$. A curve $\gamma:I\to M$, $I\subseteq\mathbb{R}$ an open interval, is a \textit{geodesic} if, and only if,
\begin{align*}
\forall\, t\in I: \nabla_{\tfrac{\mathrm{d}}{\mathrm{d}t}\gamma(t)}\tfrac{\mathrm{d}}{\mathrm{d}t}\gamma(t) = 0.
\end{align*}
A geodesic is maximal if, and only if, it does not admit a proper geodesic extension. For each $x\in M$ and $X\in\mathcal{T}_xM$, there exists a unique maximal geodesic $\gamma_X:I\to M$ with $0\in I$, $\gamma(0) = x$ and $\tfrac{\mathrm{d}}{\mathrm{d}t}\gamma_X(0) = X$. 
\end{definition-proposition}

Finally, we may recall the precise definition of the exponential map.

\begin{definition-proposition}[cf.~\cite{Lee18}]
Let $(M,g)$ be a Riemannian manifold with Levi-Civita connection $\nabla$ and geodesics $\gamma_{\cdot}$. The exponential map at $x\in M$ is defined
\begin{align*}
\exp_{x}: \set{X\in\mathcal{T}_xM\,\big\vert\,\gamma_X~\mathrm{defined~on}~[0,1]} & \to M,\\
X & \mapsto\gamma_X(1)
\end{align*}
\begin{enumerate}[(i)]
\item There exists $x\in M$ so that $\exp_x$ is defined on the whole of $\mathcal{T}_xM$ if, and only if, that is the case for all $x\in M$. In that case, $(M,g)$ is called \textit{complete}.
\item For each $x\in M$, the injectivity radius $\mathrm{inj}_g(x)$ is the largest $r\in(0,\infty]$ so that $\exp_x$ maps the ball $\mathbb{B}_g(0,r)$ with respect to $\norm{\cdot}_g$ diffeomorphically onto its image.
\item For all $r\in [0,\mathrm{inj}_g(x)]$, $\exp\big(\mathbb{B}_g(0,r)\big) = \mathbb{B}_{{d}_g}(x,r)$.
\end{enumerate}
\end{definition-proposition}

The final ingredient from Riemannian geometry used in this note is the well-known differentiability of the squared Riemannian distance.

\begin{lemma}
Let $(M,g)$ be a Riemannian manifold, and let $x\in M$. $d_g(x,\cdot)^2$ is smooth on $\mathbb{B}_{{d}_g}(x,\mathrm{inj}_g(x))$ with derivative
\begin{align}\label{eq:squared_riemannian_distance}
\mathrm{d}d_g(x,\cdot)^2 = -2\flat_{g}\exp_{(\cdot)}^{-1}(x),
\end{align}
where $\flat_{g}:\mathcal{T}M\simto\mathcal{T}^*M$ is the induced 
\end{lemma}
\begin{proof}
For $y\in\mathbb{B}_{{d}_g}(x,\mathrm{inj}_g(x))$, the definition of the exponential map gives
\begin{align*}
d_g(x,y)^2 = \langle g(x),\exp_x^{-1}(y)\otimes \exp_x^{-1}(y)\rangle.
\end{align*}
and using that the time-reversal of a geodesic is a geodesic gives
\begin{align*}
\langle\mathrm{d}_{\exp_x^{-1}(y)}\exp_x,\exp_x^{-1}(y)\rangle = -\exp_{y}^{-1}(x).
\end{align*}
Combining both formulas with the Gauss lemma, cf.~\cite[Corollary 6.2.3]{Jost17} gives~\eqref{eq:squared_riemannian_distance}.
\end{proof}

\section{The metric as a design parameter}

Let us briefly discuss the design parameter for the problem. In Euclidean funnel control, the diameter of the funnel at each time is the single design parameter, with the distance function being fixed and given by the Euclidean metric. The diameter remains a design parameter in the general case.

In Riemannian funnel control, however, an argument could be made that the metric itself becomes an additional design parameter. As enticing as it is, this is a treacherous path to tread: The metric uniquely determines the Riemannian exponential map, whose inverse takes the r\^{o}le of the difference in Euclidean funnel control. Hence, the variation of the metric determines a variation of the exponential and its inverse, which may be notoriously hard to calculate. 

On the other hand, the metric may be used to intrinsically model relaxation and strengthening of the bounds on the transient behaviour depending on the state of the reference, e.g. the overshooting of a target concentration of chemical near a critical concentration must be much more prohibited than the undershooting or the overshooting far below the critical concentration. In particular, the intrinsic modelling through the metric allows for the reference to be changed without the change of the funnel diameter.

In the following, we shall treat the Riemannian metric as part of the model; and we shall illustrate the r\^{o}le of the metric as design parameter with an example.

\section{Assumptions}

In the Euclidean funnel control literature, comparatively few assumptions are made on the error and the funnel boundary. In the first proofs, the reference trajectory is just assumed to be differentiable with locally bounded derivative, and the funnel boundary is assumed to be continuous and bounded from below by a positive constant.

For simplicity of exposition, we consider the case of smooth reference trajectories and funnel boundaries. Even then, we must impose a few assumptions which enable us to generalise the standard methods from funnel control.

First, the classical proof uses the derivative of then squared distance of the output to the reference, which is to be tracked, to prove that the latter is non-increasing in a boundary region of the funnel. On the Euclidean space, the squared distance is globally smooth; in general, the squared Riemannian distance is differentiable within the possibly finite injectivity radius. Therefore, we impose the following assumption onto the system.

\begin{assumption}\label{ass:funnel_is_small}
The reference trajectory $x_{\mathrm{ref}}\in\mathcal{C}^\infty([0,\infty),M)$ and the funnel boundary $\psi\in\mathcal{C}^\infty([0,\infty),\mathbb{R}_{>0})$ satisfy
\begin{align*}
\forall\, t\in [0,\infty): \psi(t)<\min\set{\mathrm{inj}_g\big(x_{\mathrm{ref}}(t)\big),\mathrm{dist}_g\partial_gM},
\end{align*}
where $\partial_g M := \overline{M}^g\setminus M$ and $\overline{M}^g$ denotes the completion of $M$ with respect to $d_g$.
\end{assumption}

In case of the whole Euclidean coordinate space with its Euclidean metric, the injectivity radius and the distance to the boundary are infinite, and  Assumption~\ref{ass:funnel_is_small} is automatically satisfied. Let us illustrate the non-degenerate case of the assumption with an example.

\begin{example}
Consider the Euclidean space with the origin removed, $M := \mathbb{R}^n\setminus\set{0}$, and equip $M$ with the relative Euclidean metric. The injectivity radius is the distance to the origin, so that the funnel boundary $\psi$ must satisfy $\psi(t)<\norm{x_{\mathrm{ref}}(t)}$ for all $t\in [0,\infty)$. The reason is somewhat clear: the system
\begin{align*}
\tfrac{\mathrm{d}}{\mathrm{d}t}x = c,\qquad c\neq 0,
\end{align*}
has finite escape time for all initial values in the ray $(-\infty,0)c$, so that the mere evolution of $x$ within a funnel including the origin does not \textit{a priori} guarantee infinite interval of existence.
\end{example}

A second fact that is utilised in classical funnel control is that the region that is encompassed by the funnel around the reference trajectory in finite time is compact.

\begin{assumption}\label{ass:local_boxes}
For each $T\in[0,\infty)$, the closure of the set
\begin{align}\label{eq:funnel_region}
\mathcal{F}_T := \set{(t,x)\in [0,T]\times M\,\big\vert\,d_g(x,x_{\mathrm{ref}}(t))<\psi(t)}
\end{align}
is contained in a compact set.
\end{assumption}

If the Riemannian metric $g$ is complete, then Assumption~\ref{ass:local_boxes} is automatically satisfied.

\begin{lemma}
$\mathcal{F}_T$ defined in~\eqref{eq:funnel_region} is bounded.
\end{lemma}
\begin{proof}
For any two points $(t,x),(t',x')\in\mathcal{F}_T$, we have
\begin{align*}
& d\big((t,x),(t',x')\big)\\
&\quad = \abs{t-t'}+d_g(x,x')\\
&\quad \leq T + \max_{\tau\in [0,T]}d_g\big(x_{\mathrm{ref}}(0),x_{\mathrm{ref}}(\tau)\big) + 2\max_{\tau\in [0,T]}\psi(\tau).
\end{align*}
This proves that $\mathcal{F}_T$ is bounded.
\end{proof}

For incomplete metrics, Assumption~\ref{ass:local_boxes} is not automatically satisfied as the following example demonstrates.

\begin{example}
Consider the plane with removed origin, $M = \mathbb{R}^2\setminus\set{0}$, and the induced metric. Consider further a constant reference $x_{\mathrm{ref}}\equiv c$ and a constant funnel boundary $\psi\equiv\norm{c}$. In particular, Assumption~\ref{ass:funnel_is_small} is satisfied. Then, $\mathcal{F}_T = [0,T]\times\mathbb{B}(c,\norm{c})$ is contained in a compact set if, and only if,  $\overline{\mathbb{B}}(c,\norm{c})\subseteq M$, where the closure is taken in the Riemannian closure $\overline{M} = \mathbb{R}^2$. The latter is not the case as $0\in \overline{\mathbb{B}}(c,\norm{c})$. 
\end{example}

For funnel control to not only achieve the control objective, but also have finite gain over the infinite time horizon\footnote{All real-world systems are observable on finite, but \textit{a priori} arbitrarily large, time-horizons only. For each finite time-horizon, the finiteness-assumption is not necessary to have the finite gain condition. Moreover, the size of the gain is measured in the Riemannian metric, which itself depends on the model, so that the assumption excludes, or hinders, the most interesting cases, where the reference approaches (in an ambient incomplete metric) the boundary of the complete manifold.}, so that the following finite-ness assumption may be placed on the system.

\begin{assumption}\label{ass:finite}
There exists $\varepsilon>0$ with
\begin{align*}
\forall\, t\in[0,\infty):\min\set{\psi(t),\frac{1}{\abs{\tfrac{\mathrm{d}}{\mathrm{d}t}\psi(t)}}}\geq \varepsilon,
\end{align*}
where $\frac{1}{0} := \infty$, and
\begin{align*}
\forall (t,x)\in\bigcup_{T\geq 0}\mathcal{F}_T: \max\set{\norm{\tfrac{\mathrm{d}}{\mathrm{d}t}x_{\mathrm{ref}}(t)}_g,\norm{f(x)}_g}\leq\frac{1}{\varepsilon}.
\end{align*}
\end{assumption}

\section{Funnel control on Riemannian manifolds}

In this section, we shall demonstrate that the usual funnel state feedback for the strict relative degree one vector fields gives the usual expected result, i.e. the error evolves within the funnel and is locally (globally with Assumption~\ref{ass:finite}) uniformely bounded away from the funnel boundary.

\subsection{The normal form}

By strict relative degree one, $B$ is a vector bundle isomorphism and hence there exists an input transformation which transforms $B$ into the identity, thereby giving a normal form of the control affine systems. The following theorem treats this normal form; in Proposition~\ref{prop:garbo} and Theorem~\ref{thm:funnel_high_gain}, we relax the assumption that $B^{-1}$ is known to the knowledge of a known sign-definite transformation with and without known sign.

\begin{theorem}\label{thm:funnel_normal}
Let $(M,g)$ be a Riemannian manifold, $f\in\mathfrak{X}(M)$, $x_{\mathrm{ref}}\in\mathcal{C}^\infty([0,\infty),M)$ and $\psi\in\mathcal{C}^\infty([0,\infty),\mathbb{R}_{>0})$. Let Assumptions~\ref{ass:funnel_is_small} and~\ref{ass:local_boxes} be satisfied. The closed loop
\begin{equation}\label{eq:closed_loop_I}
\begin{aligned}
\tfrac{\mathrm{d}}{\mathrm{d}t}x & = f(x)+k(\cdot)\exp_{x}^{-1}(x_{\mathrm{ref}}(\cdot))\\
k & = \frac{1}{\psi(\cdot)-d_g(x,x_{\mathrm{ref}}(\cdot))}
\end{aligned}
\end{equation}
is forward complete\footnote{A vector field $X$ is forward complete on a manifold $N$ if, and only if, each initial datum $x(t_0)\in N$ lies on an integral curve $x$ of $X$ on $[t_0,\infty)$, cf.~\cite{AngeSont99}.} on the domain
\begin{align*}
\mathcal{D} := \set{(t,x)\in [0,\infty)\times M\,\big\vert\,d_g(x,x_{\mathrm{ref}}(t))-\psi(t)<0}
\end{align*}
and each integral curve $x\in\mathcal{C}^\infty([t_0,\infty),M)$ of~\eqref{eq:closed_loop_I} satisfies
\begin{equation}\label{eq:locally_bounded_away_from_boundary}
\begin{aligned}
& \forall\, t>t_0~\exists\,\delta(T)>0~\forall\, t\in[t_0,T]: \\
& \hspace*{3cm} d_g(x(t),x_{\mathrm{ref}}(t))\leq\psi(t)-\delta(T).
\end{aligned}
\end{equation}
If additionally Assumption~\ref{ass:finite} is satisfied, then a time-invariant $\delta>0$ takes the place of $\delta(T)$ in~\eqref{eq:locally_bounded_away_from_boundary}.
\end{theorem}
\begin{proof}
The domain $\mathcal{D}$ is open in $[0,\infty)\times M$ as the preimage of a continuous function. {Moreover, due to Assumption~\ref{ass:funnel_is_small}, the righthand-side of~\eqref{eq:closed_loop_I} is smooth on $\mathcal{D}$.} Therefore, the standard existence theorem of {Picard and Lindel\"of} implies, for each initial datum $(t_0,x(t_0))\in\mathcal{D}$, the existence of a {unique} integral curve $x:[t_0,\omega)\to M$ with $(t,x(t))\in\mathcal{D}$ for all $t\in[t_0,\omega)$ and $\omega\in (t_0,\infty]$ {is maximal}. Define the error
\begin{align*}
e:[t_0,\omega)\to M,\qquad t\mapsto d_g(x(t),x_{\mathrm{ref}}(t)).
\end{align*}
Since the squared Riemannian distance is smooth on the injectivity radius of the exponential, $\frac{1}{2}e^2$ is differentiable with derivative (by chain rule)
\begin{align*}
\tfrac{\mathrm{d}}{\mathrm{d}t}\frac{1}{2}e^2(t) & = -g\big(\exp_{x(t)}^{-1}x_{\mathrm{ref}}(t),f(x(t))\big)-k(t)e(t)^2\\
& \qquad -g\big(\exp_{x(t)}^{-1}x_{\mathrm{ref}}(t),\tfrac{\mathrm{d}}{\mathrm{d}t}x_{\mathrm{ref}}(t)\big)
\end{align*}
for all $t\in[t_0,\omega)$. Let $T\in(t_0,\omega)$. By assumption~\ref{ass:local_boxes}, there exists a compact set $K(T)\subseteq \mathbb{R}\times M$ so that 
\begin{align*}
\mathrm{cl}\big(\mathcal{D}\cap [0,T]\times M\big)\subseteq K(T).
\end{align*}
Put
\begin{align*}
F(T) & := \max_{(t,x)\in K(T)}\norm{f(x)}_g,\\
R(T) & := \max_{t\in [0,T]}\norm{\tfrac{\mathrm{d}}{\mathrm{d}t}x_{\mathrm{ref}}(t)}_g.
\end{align*}
By fibrewise Cauchy-Schwarz inequality, we have
\begin{align*}
\tfrac{\mathrm{d}}{\mathrm{d}t}\frac{1}{2}e^2(t) & \leq -k(t)e(t)^2+\abs{g\big(\exp_{x(t)}^{-1}x_{\mathrm{ref}}(t),f(x(t))\big)}\\
& \qquad + \abs{g\big(\exp_{x(t)}^{-1}x_{\mathrm{ref}}(t),\tfrac{\mathrm{d}}{\mathrm{d}t}x_{\mathrm{ref}}(t)\big)}\\
& \leq -k(t)e(t)^2+{e(t)}\left(\norm{f(x(t))}_g + \norm{\tfrac{\mathrm{d}}{\mathrm{d}t}x_{\mathrm{ref}}(t)}_g\right)\\
& \leq -{k(t)}e(t)^2+{e(t)}\left(F(T) + R(T)\right).
\end{align*}
In particular, we conclude that
\begin{align*}
\tfrac{\mathrm{d}}{\mathrm{d}t}e(t)\leq F(T) + R(T) - k(t)e(t)
\end{align*}
for all $t\in[0,T]$ with $e(t)\neq 0$. Define
\begin{align}
\Psi(T) & := \abs{\min_{t\in [0,T]}\tfrac{\mathrm{d}}{\mathrm{d}t}\psi(t)},\notag\\
\psi_{\mathrm{min}}(T) & := \min_{t\in [0,T]}\psi(t),\label{eq:the_psi_min}\\
\delta_0(T) & := \min\set{1,\frac{1}{F(T)+R(T)+\Psi(T)}}\frac{\psi_{\mathrm{min}}(T)}{2}\label{eq:the_delta_zero}
\end{align}
where $\frac{1}{0} := \infty$. If $e(t)\geq\psi(t)-\delta_0(T)>0$, then
\begin{align*}
\tfrac{\mathrm{d}}{\mathrm{d}t}e(t) & \leq F(T) + R(T) - \frac{e(t)}{\psi(t)-e(t)}\\
& \leq F(T) + R(T)-\frac{\frac{\psi_{\mathrm{min}}}{2}}{\delta_0(T)} \leq -\Psi(T)\leq\tfrac{\mathrm{d}}{\mathrm{d}t}\psi(t)
\end{align*}
and therefore $e(t)\leq\psi(t)-\delta(T)$ with $\delta(T) := \min\set{\psi(t_0)-e(t_0),\delta_0(T)}$ for all $t\in[t_0,T]$. Seeking a contradiction, assume that $\omega$ is finite. Observe that then $\delta(\omega)>0$ is well-defined, and we have $e(t)\leq \psi(t)-\delta(\omega)$ for all $t\in [0,\omega)$. However, the set
\begin{small}
\begin{align*}
K(\omega) := \set{(t,x)\in [0,\omega]\times M\,\big\vert\,d_g(x,x_{\mathrm{ref}}(t))\leq\psi(t)-\delta(\omega)}
\end{align*}
\end{small}
\noindent is a compact subset of $\mathcal{D}$ so that, since $x$ is {forward} maximal, $x(t_1)\notin K(\omega)$ for some $t_1\in[t_0,\omega)$, contradicting $x(t)\in K(\omega)$ for all $t\in[t_0,\omega)$. Therefore, $\omega = \infty$ and~\eqref{eq:locally_bounded_away_from_boundary} is satisfied. {If Assumption~\ref{ass:finite} is satisfied, then the time-dependent constants $F(T)$, $R(T)$, $\psi(T)$ and $\psi_{\mathrm{min}}(T)$ may be replaced by time-invariant constants and hence $\delta$ is independent of time.}
\end{proof}

\subsection{The general form}

As usual in funnel control, the input field $B$ is, in general un-known, so that Theorem~\ref{thm:funnel_normal} is rather restrictive. In case that the sign of the gain function is known, we {sketch an adapted} proof.

\begin{proposition}\label{prop:garbo}
Let $(M,g)$ be a Riemannian manifold, let $f\in\mathfrak{X}(M)$ and let $B\in\Gamma(\mathcal{T}M\otimes \mathcal{T}^*M)$ be a vector bundle automorphism so that
\begin{align*}
\forall\, x\in M\ \exists\, \gamma>0\ \forall\, X\in \mathcal{T}_xM:\ g(X,\langle B,X\rangle)\geq \gamma g(X,X).
\end{align*}
Let Assumptions~\ref{ass:funnel_is_small} and~\ref{ass:local_boxes} be satisfied. The closed loop
\begin{equation}\label{eq:closed_loop_II}
\begin{aligned}
\tfrac{\mathrm{d}}{\mathrm{d}t}x & = f(x)+k(\cdot)B(x)\exp_{x}^{-1}(x_{\mathrm{ref}}(\cdot))\\
k & = \frac{1}{\psi(\cdot)-d_g(x,x_{\mathrm{ref}}(\cdot))}
\end{aligned}
\end{equation}
{is forward complete on the domain}
\begin{align*}
\mathcal{D} := \set{(t,x)\in [0,\infty)\times M\,\big\vert\,d_g(x,x_{\mathrm{ref}}(t))-\psi(t)<0}
\end{align*}
and {each integral curve $x\in\mathcal{C}^\infty([t_0,\infty),M)$ of~\eqref{eq:closed_loop_II}} satisfies~\eqref{eq:locally_bounded_away_from_boundary}. If additionally Assumption~\ref{ass:finite} is satisfied, then a time-invariant $\delta>0$ takes the place of $\delta(T)$ in~\eqref{eq:locally_bounded_away_from_boundary}.
\end{proposition}
\begin{proof}
In the following, we shall make use of the notation and objects introduced in the proof of Theorem~\ref{thm:funnel_normal}. Consider the Riemannian metric
\begin{align*}
\widehat{g}:\mathcal{T}M\otimes\mathcal{T}M & \to\mathbb{R},\\
\widehat{g}(X,Y) & := \frac{1}{2}\big(g(X,\langle B,Y\rangle)+g(Y,\langle B,X\rangle)\big).
\end{align*}
Then, the derivative of the squared Riemannian distance $\frac{1}{2}e^2(\cdot) = \frac{1}{2}d_g(x(\cdot),x_{\mathrm{ref}}(\cdot))^2$ reads as
\begin{align*}
\tfrac{\mathrm{d}}{\mathrm{d}t}\frac{1}{2}e^2(t) & = -g\big(\exp_{x(t)}^{-1}x_{\mathrm{ref}}(t),f(x(t))\big)\\
& \qquad -k(t)\norm{\exp_{x(t)}^{-1}x_{\mathrm{ref}}(t)}_{\widehat{g}}^2\\
& \qquad -g\big(\exp_{x(t)}^{-1}x_{\mathrm{ref}}(t),\tfrac{\mathrm{d}}{\mathrm{d}t}x_{\mathrm{ref}}(t)\big),
\end{align*}
and $\norm{\exp_{x(t)}^{-1}x_{\mathrm{ref}}(t)}_{\widehat{g}}^2\geq \gamma(T) e(t)^2$ with
\begin{align*}
\gamma(T) := \inf_{(x,t)\in K(T)\setminus\mathrm{graph}(x_{\mathrm{ref}})}\frac{\norm{\exp_{x_{\mathrm{ref}(t)}}^{-1}(x)}_{\widehat{g}}}{\norm{\exp_{x_{\mathrm{ref}(t)}}^{-1}(x)}_g}>0
\end{align*}
gives the inequality
\begin{align*}
\tfrac{\mathrm{d}}{\mathrm{d}t}e(t)\leq F(T)+R(T)-\gamma(T) k(t)e(t).
\end{align*}
Multiplying $\delta_0$ with $\frac{1}{\gamma(T)}$ and proceeding as in the proof of Theorem~\ref{thm:funnel_normal} gives the claimed result.
\end{proof}

{The most general case cannot be treated with the funnel controller. It must} be assumed that the gain matrix is sign-definite, i.e. $B$ is everywhere positive or negative definite.

\begin{assumption}\label{ass:sign_definite}
There exists a known vector bundle isomorphism $\digamma:\mathcal{T}M\simto E$ so that $B\digamma$ satisfies
\begin{align}\label{eq:sign_definite}
\forall X\in\mathcal{T}_xM\setminus\set{0}: \sigma g(X,B\digamma X)>0 
\end{align}
for all $x$ in an open neighbourhood of the funnel region $\set{x\in M\,\big\vert\,\exists\, t\in[0,\infty): d_g(x,x_{\mathrm{ref}}(t))<\psi(t)}$ and $\sigma\in\set{\pm 1}$ which is independent of $x$.
\end{assumption}

\begin{remark}
It is not necessary to \textit{a priori} claim that $\sigma$ in Assumption~\ref{ass:sign_definite} is independent of $x$. Indeed, in any local coordinates defined on an open and connected set $U$, $g\vert_U$ corresponds to a smooth pointwise symmetric and positive definite matrix field $g\vert_U\in\mathcal{C}^\infty(U,\mathbb{R}^{n\times n})$, and $(B\digamma)\vert_U$ corresponds to a smooth matrix field $(B\digamma)\vert_U\in\mathcal{C}^\infty(U,\mathbb{R}^{n\times n})$. The assumption is then that, for each $x\in U$, $g(x)B(x)\digamma(x)+\digamma(x)^\top B(x)^\top g(x)$ is positive definite or negative definite. Since the sets of positive and negative definite (symmetric) matrices are open and disjoint, and since $g(\cdot)B(\cdot)\digamma(\cdot)+\digamma(\cdot)^\top B(\cdot)^\top g(\cdot)$ is, in particular, continuous, $g(x)B(x)\digamma(x)+\digamma(x)^\top B(x)^\top g(x)$ is either positive definite ($\sigma = 1$) for all $x\in U$ or negative definite ($\sigma = -1$) for all $x\in U$. Since $M$ is connected (as we assume throughout the note), we conclude \textit{a posteriori} that $\sigma$ is indeed independent of $x$.
\end{remark}

\begin{theorem}\label{thm:funnel_high_gain}
Let $(M,g)$ be a Riemannian manifold, let $f\in\mathfrak{X}(M)$ and let $N\in\mathcal{C}^\infty([0,\infty),\mathbb{R})$ be surjective. Let Assumptions~\ref{ass:funnel_is_small} and~\ref{ass:local_boxes} be satisfied, and let $\digamma\in\Gamma(E\otimes \mathcal{T}^*M)$ be a vector bundle isomorphism realising assumption~\ref{ass:sign_definite}. The closed loop
\begin{equation}\label{eq:closed_loop_III}
\begin{aligned}
\tfrac{\mathrm{d}}{\mathrm{d}t}x & = f(x)+N(k(\cdot))B(x)\digamma(x)\exp_{x}^{-1}(x_{\mathrm{ref}}(\cdot))\\
k & = \frac{1}{\psi(\cdot)-d_g(x,x_{\mathrm{ref}}(\cdot))}
\end{aligned}
\end{equation}
{is forward complete} on the domain
\begin{align*}
\mathcal{D} := \set{(t,x)\in [0,\infty)\times M\,\big\vert\,d_g(x,x_{\mathrm{ref}}(t))-\psi(t)<0}
\end{align*}
and {each integral curve $x\in\mathcal{C}^\infty([t_0,\infty),M)$ of~\eqref{eq:closed_loop_III}} satisfies~\eqref{eq:locally_bounded_away_from_boundary}. If additionally Assumption~\ref{ass:finite} is satisfied and $B\digamma$ is uniformely bounded away from zero\footnote{i.e. the zero in the inequality~\eqref{eq:sign_definite} is replaced by a positive constant independent of $x$.}, then a time-invariant $\delta>0$ takes the place of $\delta(T)$ in~\eqref{eq:locally_bounded_away_from_boundary}.
\end{theorem}
\begin{proof}
Let $x:[t_0,\omega)\to M$ be a{n integral curve} of~\eqref{eq:closed_loop_III} with {maximal} $\omega\in(t_0,\infty]$, and put $e(\cdot) := d_g(x(\cdot),x_{\mathrm{ref}}(\cdot))$. Let $T\in (t_0,\omega)$, put $\widehat{g}:= B\digamma$, and consider $K(T)$, $F(T)$ and $R(T)$ as in the proof of Theorem~\ref{thm:funnel_normal}. Without loss of generality, $K(T)$ may be chosen so that~\eqref{eq:sign_definite} holds on the whole of {the Cartesian projection} $p_M(K(T))$. Then, by Assumption~\ref{ass:sign_definite}, there is $\gamma(T)>0$ so that 
\begin{align*}
\forall (t,x)\in K(T)\forall X\in\mathcal{T}_xM: \sigma g(X,\widehat{g}(x)X)\geq\gamma(T) g(X,X).
\end{align*}
and hence, we have the inequality
\begin{align*}
\tfrac{\mathrm{d}}{\mathrm{d}t}\frac{1}{2}e^2(t) & \leq F(T)+R(T) - \sigma N(k(t)) \frac{\sigma}{e(t)} \widehat{e}(t)
\end{align*}
for all $t\in [t_0,T]$ with $e(t)\neq 0$, where
\begin{align*}
\widehat{e}(t) := g\big(\exp_{x(t)}^{-1}x_{\mathrm{ref}}(t),\widehat{g}(x(t))\exp_{x(t)}^{-1}x_{\mathrm{ref}}(t)\big).
\end{align*}
In particular, we have $\frac{\sigma}{e(t)} \widehat{e}(t)\geq\gamma(T)e(t)^2$. Moreover, since $N$ is surjective, there are $\max\set{k(0),\frac{2}{\psi_{\mathrm{min}}(T)}}< z_1<z_2$ with
\begin{align*}
\forall z\in [z_1,z_2]:  \sigma N(z)> \frac{1}{\gamma(T)\delta_0(T)}
\end{align*}
where $\psi_{\mathrm{min}}(T)$ and $\delta_0(T)$ are defined as in~\eqref{eq:the_psi_min} and~\eqref{eq:the_delta_zero}. Seeking a contradiction, assume that $k(t)> z_1$ for some $t\in [0,T]$; in particular, $e(t)>0$. Then, we conclude analogously as in the proof of Theorem~\ref{thm:funnel_normal} that
\begin{align*}
\tfrac{\mathrm{d}}{\mathrm{d}t}e(t) & \leq F(T)+R(T)-\frac{1}{\gamma(T)\delta_0(T)}\gamma(T)e(t)\\
& < F(T)+R(T)-\frac{1}{\gamma(T)\delta_0(T)}\gamma(T)\frac{\psi_\mathrm{min}(T)}{2}\\
& \leq \tfrac{\mathrm{d}}{\mathrm{d}t}\psi(t)
\end{align*}
and hence $\tfrac{\mathrm{d}}{\mathrm{d}t}k(t)<0$. This shows that $k$ is bounded on $[t_0,T]$, and analogously as to the proof of Theorem~\ref{thm:funnel_normal}, we conclude that~\eqref{eq:locally_bounded_away_from_boundary} holds and $\omega = \infty$.
\end{proof}

\section{Example: Slow dynamics of a pendulum}

Consider a massless rod of unit length attached at one end to a motor fixed in space, and consider a mass $m>0$ attached at the free end. The position of the mass may be modelled as a point $[\vartheta]\in S^1 = \mathbb{R}/2\pi\mathbb{Z}$ represented by values in $[-\pi,\pi]$ with $\pi$ and $-\pi$ identified. With drag coefficient $c>0$, the gravity constant $g$, we get 
\begin{align*}
\tfrac{\mathrm{d}^2}{\mathrm{d}t^2}\vartheta = -g\sin\vartheta-\frac{c}{m}\cos\tfrac{\mathrm{d}}{\mathrm{d}t}\vartheta + \upsilon,
\end{align*}
where $u$ is the force acting on the mass by the motor. When $\frac{c}{m}$ is sufficiently large (i.e. $m$ is sufficiently small and $c$ is sufficiently large, e.g. a bushy piece of cotton or a down feather), the evolution may be approximated by the slow evolution
\begin{align}\label{eq:slow_dynamics}
\tfrac{\mathrm{d}}{\mathrm{d}t}\vartheta = -\frac{mg}{c}\sin\vartheta + u
\end{align}
with $u:=\frac{1}{c}\upsilon$. Then,~\eqref{eq:slow_dynamics} is treated by  Theorem~\ref{thm:funnel_normal}. Consider the reference trajectory
\begin{align*}
\vartheta_{\mathrm{ref}}:[0,\infty)\to S^1,\qquad t\mapsto \frac{\pi}{2}\sin(t).
\end{align*}
Recall that a Riemannian metric is given by an  everywhere positive $2\pi$-periodic function $g\in\mathcal{C}^\infty(\mathbb{R},\mathbb{R}_{>0})$. Put $2L = \int_{-\pi}^{\pi}g(\vartheta)\,\mathrm{d}\vartheta$, so that the injectivity radius at each point is $L$; the boundary of $S^1$ is empty. Then, the admissible funnel boundaries $\psi$ satisfying Assumption~\ref{ass:funnel_is_small} are characterised by $\psi<L$, and Assumption~\ref{ass:local_boxes} is automatically satisfied as $S^1$ is compact. For example, consider the metric given by
\begin{align*}
g(\vartheta) := 0.1 + \frac{1 + \cos(\vartheta)}{2\pi}.
\end{align*}
This metric is symmetric around zero, where it has its global maximum, but not symmetric elsewhere. The global minimum is at $\pi = -\pi$. Consider the constant funnel boundary $\psi\equiv 0.1$. Then, the metric enforces stricter tracking near zero and looser tracking near $\pi$. Moreover, the funnel is opened more away from zero, so overshooting is punished less than lagging behind. {Observe that the stricter and looser tracking of certain reference states is determined by the minimum and maximum of the metric, which may be tuned to achieve the desired behaviour.}

\begin{figure}[h]
\centering
\includegraphics[scale=0.5]{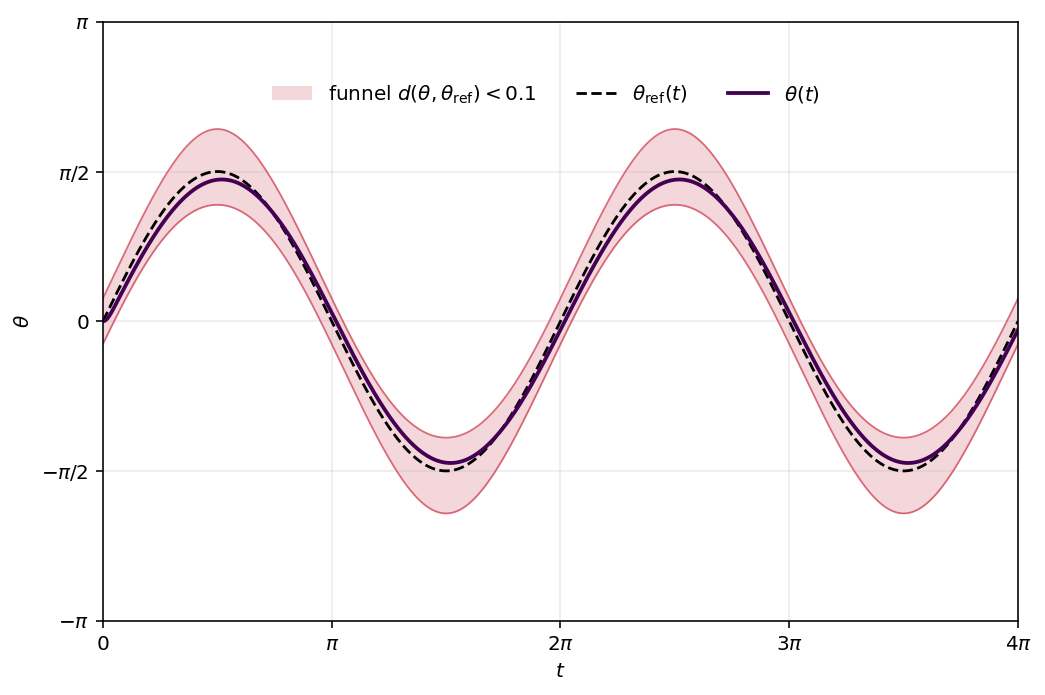}
\caption{Reference tracking}
\label{fig:garbage}
\end{figure}

Given $\vartheta\in [-\pi,\pi]$ and $x\in\mathbb{R}$, the geodesic $\gamma_x$ in $\vartheta$ satisfies
\begin{align*}
\int_{\vartheta}^{\gamma(t)}g(\alpha)\,\mathrm{d}\alpha = tx g(\vartheta)
\end{align*}
and hence, for $\theta'$ with $d_g(\theta,\theta')<L$, we obtain
\begin{align*}
\exp_{\vartheta}^{-1}\vartheta' = \frac{1}{g(\vartheta)}\int_{\vartheta}^{\vartheta'}g(\alpha)\,\mathrm{d}\alpha.
\end{align*}
Therefore, the funnel controller may be implemented directly. The result is presented in Figure~\ref{fig:garbage}.

\section{Conclusion and outlook}

In this note, we proposed an extension of the classical Euclidean funnel controller to Riemannian manifolds. For simplicity of exposition, we studied control-affine systems with relative degree one, which takes the r\^{o}le of the linear prototype of Euclidean funnel control.

We generalised the Euclidean funnel controller to the Riemannian funnel controller
\begin{align*}
u(t) = N\left(\frac{1}{\psi(t)-d_g(x(t),x_{\mathrm{ref}}(t))}\right)\exp_x^{-1}\big(x_{\mathrm{ref}}(t)\big)
\end{align*}
and adapted the classical proof to show that under the assumptions that guarantee that $u$ is smooth within the funnel the control objective is achieved. The Riemannian funnel controller is coordinate-free and intrinsically defined, while the Euclidean funnel controller (extended to open sets) requires the choice of normal coordinates. 

A consequence of Riemannian funnel control is that the Euclidean metric in Euclidean funnel control can be replaced by any Hilbert space norm on $\mathbb{R}^n$, since these norms are geodesically equivalent, and the control objective is achieved.

In future work, the approach will be generalised to systems with nonlinear output and arbitrary relative degree; additionally, the reference trajectory will be generalised to a fibre bundle of reference submanifolds. 

\bibliographystyle{plain}
\bibliography{Literatur}

\end{document}